\documentclass[twoside,a4paper,reqno,11pt]{amsart} 
\usepackage{amsfonts, amsbsy, amsmath, amssymb, latexsym}
\usepackage{mathrsfs,array}
\usepackage[top=30mm,right=30mm,bottom=30mm,left=30mm]{geometry}
\usepackage{stmaryrd}
\usepackage{bm}
\usepackage{xcolor}
\usepackage{hyperref}
\usepackage{listings}

\allowdisplaybreaks

\usepackage{booktabs}

\newtheorem{theorem}{Theorem}[section]

\newtheorem{proposition}[theorem]{Proposition}

{\theoremstyle{definition}

}

{\theoremstyle{definition}

}

\usepackage{mathtools}

\begin{document}
	
\author{Melissa Lee and Tomasz Popiel}
\address{Department of Mathematics, University of Auckland, Auckland, New Zealand}
\email{melissa.lee@auckland.ac.nz, tomasz.popiel@auckland.ac.nz}

\title{$\text{M}$, $\text{B}$ and $\text{Co}_1$ are recognisable by their prime graphs}

\begin{abstract} 
The {\em prime graph}, or {\em Gruenberg--Kegel graph}, of a finite group $G$ is the graph $\Gamma(G)$ whose vertices are the prime divisors of $|G|$, and whose edges are the pairs $\{p,q\}$ for which $G$ contains an element of order $pq$. 
A finite group $G$ is {\em recognisable} by its prime graph if every finite group $H$ with $\Gamma(H)=\Gamma(G)$ is isomorphic to $G$. 
By a result of Cameron and Maslova, every such group must be almost simple, so one natural case to investigate is that in which $G$ is one of the $26$ sporadic simple groups. 
Existing work of various authors answers the question of recognisability by prime graph for all but three of these groups, namely the Monster, $\text{M}$, the Baby Monster, $\text{B}$, and the first Conway group, $\text{Co}_1$.
We prove that these three groups are recognisable by their prime graphs. 
\end{abstract}

\date{\today}
\maketitle

\section{Introduction}	 \label{s:intro}

The {\em Gruenberg--Kegel graph} of a finite group $G$, introduced in an unpublished 1975 manuscript of Karl~W.~Gruenberg and Otto~H.~Kegel, is the labelled graph $\Gamma(G)$ whose vertices are the prime divisors of $|G|$, and whose edges are the pairs $\{p,q\}$ for which $G$ contains an element of order $pq$. 
It is also called the {\em prime graph} of $G$, and we use this name here for brevity. 
A finite group $G$ is {\em recognisable} by its prime graph if every finite group $H$ with $\Gamma(H)=\Gamma(G)$ is isomorphic to $G$. 
More generally, $G$ is $k$-{\em recognisable} by its prime graph if there are precisely $k$ pairwise non-isomorphic groups with the same prime graph as $G$. 
If no such $k$ exists, then $G$ is {\em unrecognisable} by its prime graph. 

The question of recognisability of various groups by their prime graphs has attracted significant interest. 
We refer the reader to the recent article of Cameron and Maslova~\cite{CamMas} for an up-to-date review of the literature (and several new results). 
Much work has focused on simple groups, and this is justified by \cite[Theorem~1.2]{CamMas}, which states, in particular, that if $G$ is $k$-recognisable by its prime graph for some $k$, then $G$ is almost simple, i.e. $G_0 \leqslant G \leqslant \text{Aut}(G_0)$ for some non-abelian simple group $G_0$. 
It is therefore natural to ask, in particular, which of the $26$ sporadic simple groups are recognisable by their prime graphs. 
As summarised in Table~\ref{tabSummary}, this question has been answered for all but three of these groups: the Monster, $\text{M}$, the Baby Monster, $\text{B}$, and the first Conway group, $\text{Co}_1$. 
The purpose of this note is to settle these remaining cases. 
We prove the following theorem. 

\begin{theorem} \label{thm1}
$\textnormal{M}$, $\textnormal{B}$ and $\textnormal{Co}_1$ are recognisable by their prime graphs. 
\end{theorem}

\begin{table}[!ht]
\small
\begin{tabular}{llll}
\toprule
$G$ & Recognisable? & $\Gamma(H)=\Gamma(G)$? & Reference \\
\midrule
$\text{Co}_2$, $\text{J}_1$, $\text{M}_{22}$, $\text{M}_{23}$, $\text{M}_{24}$ & Yes && Hagie~\cite{Hagie} \\
$\text{M}_{11}$ & $2$-recognisable & $H=\text{L}_2(11)$ & Hagie~\cite{Hagie} \\
$\text{M}_{12}$ & Unrecognisable && Hagie~\cite{Hagie} \\
$\text{J}_2$ & Unrecognisable && Shi et al.~\cite{MazurovShi,PraegerShi} \\
$\text{J}_4$ & Yes && Zavarnitsine~\cite{Zavar} \\
$\text{Ru}$ & Yes && Kondrat'ev~\cite{Kon1} \\
$\text{HN}$ & $2$-recognisable & $H = \text{HN}.2$ & Kondrat'ev~\cite{Kon1} \\
$\text{Fi}_{22}$ & $3$-recognisable & $H \in \{ \text{Fi}_{22}.2, \text{Suz}.2\}$ & Kondrat'ev~\cite{Kon1} \\
$\text{He}$, $\text{McL}$, $\text{Co}_3$ & Unrecognisable && Kondrat'ev~\cite{Kon1} \\
$\text{Fi}_{23}$, $\text{Fi}_{24}'$, $\text{J}_3$, $\text{Ly}$, $\text{O'N}$, $\text{Suz}$, $\text{Th}$ & Yes && Kondrat'ev~\cite{Kon2} \\
$\text{HS}$ & $2$-recognisable & $H=\text{U}_6(2)$ & Kondrat'ev~\cite{Kon2} \\
$\text{M}$, $\text{B}$, $\text{Co}_1$ & Yes && Theorem~\ref{thm1} \\
\bottomrule
\end{tabular}
\caption{Recognisability of the sporadic simple groups by their prime graphs.}
\label{tabSummary}
\end{table}

\section{Proof of Theorem~\ref{thm1} --- $\text{Co}_1$} \label{s:proofCo1}

Recall that $|\text{Co}_1| = 2^{21} \cdot 3^9 \cdot 5^4 \cdot 7^2 \cdot 11 \cdot 13 \cdot 23$. 
The vertex set of $\Gamma(\text{Co}_1)$ is therefore $\{2,3,5,7,11,13,23\}$, and we see from the ATLAS~\cite{ATLAS} that $\Gamma(\text{Co}_1)$ has two connected components, one of which is the isolated vertex $23$. 
We also recall from the ATLAS~\cite{ATLAS} that $\text{Co}_1$ contains both $\text{Co}_3$ and $(\text{A}_7 \times \text{L}_2(7)):2$ as maximal subgroups. 
(Here we follow the ATLAS~\cite{ATLAS} and use Artin's notation ``L'' for ``PSL'', e.g. $\text{L}_2(7) = \text{PSL}_2(7)$.)

Our argument begins along the lines of Kondrat'ev's proofs \cite{Kon1,Kon2}, which in turn rely on earlier work of Hagie~\cite{Hagie}. 
Suppose that $G$ is a finite group with $\Gamma(G) = \Gamma(\text{Co}_1)$. 
Then \cite[Theorem~3]{Hagie} implies that $G/F(G) \cong \text{Co}_1$, and that the prime divisors of the Fitting subgroup $F(G)$ of $G$ form a (possibly empty) subset of $\{2,3\}$. 
We must show that $F(G)$ is trivial, so we suppose towards a contradiction that the maximal normal $p$-subgroup $O_p(G)$ of $G$ is non-trivial for at least one $p \in \{2,3\}$. 
This implies that every $g \in G$ of order $23$ must act fixed-point freely (by conjugation) on $O_p(G)$, because otherwise the product of $g$ with a commuting $h \in O_p(G)$ would have order $23p$, so $\Gamma(G)$ would contain the edge $\{23,p\}$, in contradiction with the assumption that $\Gamma(G)=\Gamma(\text{Co}_1)$. 
Following Kondrat'ev \cite{Kon1,Kon2}, the idea is to obtain a contradiction by considering $p$-modular representations of $\text{Co}_1$ on which elements of order $23$ might act fixed-point freely, as we now explain. 

If $O_p(G)$ is non-trivial for a given $p \in \{2,3\}$, then, by suitably choosing two neighbouring terms $K$ and $L$ of a chief series for $G$, with $K < L \leqslant F(G)$, we can obtain a chief factor $V = K/L$ such that $V$ is an elementary abelian $p$-group. 
(Here $V$ is necessarily elementary abelian, being solvable and characteristically simple, and we can choose $K$ and $L$ such that $V$ is a $p$-group for the given $p$.)
Since $V$ is a minimal normal subgroup of $G/K$, we may regard it as an irreducible module for $G/K$ in characteristic $p$. 
Furthermore, since the centraliser $C_{G/K}(V)$ of $V$ in $G/K$ is equal to $F(G)/K$, the N/C theorem implies that we may regard $V$ as a {\em faithful} irreducible module for $(G/K)/(F(G)/K) \cong G/F(G) \cong \text{Co}_1$. 
(To verify that $C_{G/K}(V) = F(G)/K$, first observe that $F(G)/K = F(G/K)$. 
Next, note that $F(G/K) \leqslant C_{G/K}(V)$, because $F(G/K)$ centralises every minimal normal subgroup of $G/K$. 
Finally, note that if the inclusion was strict, then the kernel of the action of $G/K$ on $V$ would sit between $F(G/K)$ and $C_{G/K}(V)$, so the correspondence theorem would imply the existence of a non-trivial proper normal subgroup of $G/F(G) \cong \text{Co}_1$.)

The upshot of the above argument is this: under the assumption that $O_p(G)$ is non-trivial for a given $p \in \{2,3\}$, it must be the case that $\text{Co}_1 \cong G/F(G)$ admits a faithful irreducible module in characteristic $p$ on which elements of order $23$ act fixed-point freely. 
We now argue that $\text{Co}_1$ admits no such module, for either $p \in \{2,3\}$. 
In fact, we prove the following stronger result, which we record explicitly in case it is of independent interest. 

\begin{proposition}
Let $V$ be a faithful irreducible module for $\textnormal{Co}_1$ in any characteristic other than $23$. 
Then every element of order $23$ in $\textnormal{Co}_1$ fixes some point of $V$. 
\end{proposition}

First consider $p=3$. 
An element of order $23$ in $\text{Co}_1$ normalises a subgroup of order $2^{11}$ (see \cite{ATLAS}), so it follows from \cite[Theorem 3.4.4]{Gorenstein} that such an element has a fixed point on every faithful $\mathbb{F}\text{Co}_1$-module, for every field $\mathbb{F}$ of characteristic not equal to $2$ or $23$. 
In particular, such an element does {\em not} act fixed-point freely on any faithful irreducible $\mathbb{F}\text{Co}_1$-module defined over a field $\mathbb{F}$ of characteristic $3$. 
Therefore, $O_3(G)$ must be trivial. 

Now consider $p=2$. 
Suppose that $V$ is a faithful irreducible module for $\text{Co}_1$ in characteristic~$2$, and let $\chi$ be the corresponding Brauer character. 
Consider a maximal subgroup $\text{Co}_3$ of $\text{Co}_1$. 
According to the $2$-modular character table of $\text{Co}_3$, which is available in {\sf GAP}~\cite{GAPbc,GAP}, the only irreducible module for $\text{Co}_3$ in characteristic $2$ on which elements of order $23$ act fixed-point freely is the module, call it $W$, of dimension $22$. 
Therefore, each of the $k$ composition factors of the restriction $V \downarrow \text{Co}_3$ of $V$ to $\text{Co}_3$ is isomorphic to $W$, and so the restriction of $\chi$ to $\text{Co}_3$ is $k$ times the Brauer character of $W$. 
Now consider a maximal subgroup $H=(\text{A}_7 \times \text{L}_2(7)):2$ of $\text{Co}_1$. 
Note that there are exactly $8$ conjugacy classes of elements of odd order in $\text{Co}_1$ that intersect both $\text{Co}_3$ and $H$. 
These classes are listed in Table~\ref{tab1}, together with the corresponding values of the hypothesised $2$-modular Brauer character $\chi$ and the corresponding classes in $\text{Co}_3$ and $H$. 
(Note that the class fusions can be readily determined in {\sf GAP}~\cite{GAPbc,GAP} by applying the functions \texttt{FusionConjugacyClasses} and \texttt{ClassNames} to the relevant ordinary character tables.)

Now, consider that the restriction of $\chi$ to $H$ must be a linear combination of the irreducible $2$-modular Brauer characters for $H$, with non-negative integer coefficients. 
Using the $2$-modular character table of $H$, which is also available in {\sf GAP}~\cite{GAPbc,GAP}, and comparing with the final column of Table~\ref{tab1}, we obtain a system of equations for the multiplicities of the characters appearing in this restriction. 
There are $8$ equations in $16$ unknowns: one equation for each of the $8$ conjugacy classes of $H$ appearing in Table~\ref{tab1}, and one variable for each of the irreducible $2$-modular Brauer characters for $H$. 
Explicitly, if we label these variables $x_1,\dots,x_{16}$ according to the numbering of the characters for $H$ used in {\sf GAP}~\cite{GAPbc,GAP}, then we obtain the linear system with augmented matrix
\[
\small
\left[
\begin{array}{rrrrrrrrrrrrrrrr}
1&8&6&14&20&6&24&24&36&84&120&8&64&48&112&160 \\
1&2&\cdot&-1&-1&6&6&6&\cdot&-6&-6&8&16&\cdot&-8&-8 \\
1&8&6&14&20&\cdot&\cdot&\cdot&\cdot&\cdot&\cdot&-1&-8&-6&-14&-20 \\
1&-4&3&2&-4&\cdot&\cdot&\cdot&\cdot&\cdot&\cdot&-1&4&-3&-2&4 \\
1&2&\cdot&-1&-1&\cdot&\cdot&\cdot&\cdot&\cdot&\cdot&-1&-2&\cdot&1&1 \\
1&1&-1&\cdot&-1&6&3&3&-6&\cdot&-6&8&8&-8&\cdot&-8 \\
1&1&-1&\cdot&-1&-1&-4&3&1&\cdot&1&1&1&-1&\cdot&-1 \\
1&1&-1&\cdot&-1&\cdot&\cdot&\cdot&\cdot&\cdot&\cdot&-1&-1&1&\cdot&1
\end{array}
\right| \hspace{-2pt}
\left.
\begin{array}{r}
22k \\ 4k \\ -2k \\ -2k \\ -2k \\ k \\ k \\ -2k
\end{array} 
\right],
\]
where $\cdot$ denotes $0$. 
As explained below, this system has no non-negative integer solutions $(x_1,\dots,x_{16})$ for any positive integer $k$, so the hypothesised module $V$ for $\text{Co}_1$ does not exist. 
Therefore, $O_2(G)$ is also trivial, so $F(G)$ is trivial and hence $G \cong \text{Co}_1$, as claimed. 

It remains to verify that the aforementioned linear system has no non-negative integer solutions. 
It suffices to consider $k=1$ and check that there are no non-negative rational solutions. 
When $k=1$, the augmented matrix given above is row equivalent to
\[
\left[
\begin{array}{rrrrrrrrrrrrrrrr}
1&\cdot&\cdot&\cdot&-2&\cdot&\cdot&\cdot&\cdot&\cdot&\cdot&-1&\cdot&\cdot&\cdot&2 \\
\cdot&1&\cdot&\cdot&1&\cdot&\cdot&\cdot&\cdot&\cdot&\cdot&\cdot&-1&\cdot&\cdot&-1 \\
\cdot&\cdot&1&\cdot&\cdot&\cdot&\cdot&\cdot&\cdot&\cdot&\cdot&\cdot&\cdot&-1&\cdot&\cdot \\
\cdot&\cdot&\cdot&1&1&\cdot&\cdot&\cdot&\cdot&\cdot&\cdot&\cdot&\cdot&\cdot&-1&-1 \\
\cdot&\cdot&\cdot&\cdot&\cdot&1&\cdot&\cdot&-2&\cdot&-2&3/2&\cdot&-3&\cdot&-3 \\
\cdot&\cdot&\cdot&\cdot&\cdot&\cdot&1&\cdot&1&\cdot&1&-1/2&1&2&\cdot&2 \\
\cdot&\cdot&\cdot&\cdot&\cdot&\cdot&\cdot&1&1&\cdot&1&1/2&2&1&\cdot&1 \\
\cdot&\cdot&\cdot&\cdot&\cdot&\cdot&\cdot&\cdot&\cdot&1&1&\cdot&\cdot&\cdot&3/2&3/2
\end{array}
\right| \hspace{-2pt}
\left.
\begin{array}{r}
-2 \\ \cdot \\ \cdot \\ \cdot \\ \cdot \\ \cdot \\ 1 \\ \cdot
\end{array} 
\right].
\]
Given that all of $x_1,\dots,x_{16}$ should be non-negative, the final row immediately implies that $x_{10}=x_{11}=x_{15}=x_{16}=0$, and then the fourth row yields $x_4=x_5=0$. 
The first row then implies that $x_{12} = x_1+2$, so in particular $x_{12} \geqslant 2$. 
Putting $x_{12} \geqslant 2$ and $x_{11}=x_{16}=0$ into the second-last row then forces $x_8+x_9+2x_{13}+x_{14} \leqslant 0$, so all of $x_8$, $x_9$, $x_{13}$ and $x_{14}$ must also vanish. 
Looking again at the second-last row, it follows that $x_{12}=2$. 
However, the fifth row then implies that $x_6 = -3 < 0$. 

\begin{table}[!t]
\begin{tabular}{lllr}
\toprule
$\text{Co}_1$ & $\text{Co}_3$ & $(\text{A}_7 \times \text{L}_2(7)):2$ & $\chi$ \\
\midrule
1A & 1A & 1A & $22k$ \\
3B & 3B & 3D & $4k$ \\
3D & 3C & 3A, 3C, 3E & $-2k$ \\
7B & 7A & 7B, 7D & $k$ \\
21C & 21A & 21C & $-2k$ \\
\bottomrule
\end{tabular}
\caption{The conjugacy classes of odd-order elements in $\text{Co}_1$ that intersect both of the maximal subgroups $\text{Co}_3$ and $H=(\text{A}_7 \times \text{L}_2(7)):2$, and the corresponding conjugacy classes in these subgroups. 
Class names are as in the ATLAS~\cite{ATLAS} and {\sf GAP}~\cite{GAPbc,GAP}. 
The final column lists the corresponding values of the $2$-modular Brauer character $\chi$ of the hypothesised module $V$ for $\text{Co}_1$ considered in Section~\ref{s:proofCo1}, which supposedly restricts to $\text{Co}_3$ as $k$ copies of the unique $22$-dimensional module for $\text{Co}_3$, for some $k$.}
\label{tab1}
\end{table}

\section{Proof of Theorem~\ref{thm1} --- $\text{M}$} \label{s:proofM}

Recall that $|\text{M}| = 2^{46} \cdot 3^{20} \cdot 5^9 \cdot 7^6 \cdot 11^2 \cdot 13^3 \cdot 17 \cdot 19 \cdot 23 \cdot 29 \cdot 31 \cdot 41 \cdot 47 \cdot 59 \cdot 71$. 
We see from the ATLAS~\cite{ATLAS} that the prime graph $\Gamma(\text{M})$ has three isolated vertices, namely $41$, $59$ and $71$, and that all of its other vertices comprise a connected component. 
We also recall from the ATLAS~\cite{ATLAS} that $\text{M}$ contains a maximal subgroup of the form $3^8.\text{O}_8^-(3).2_3$ (a non-split extension of the orthogonal group $\text{O}_8^-(3).2_3$ acting naturally on $\mathbb{F}_3^8$, where the $2_3$ refers to one of three possible ways of extending $\text{O}_8^-(3)$ to a group of the form $\text{O}_8^-(3).2$). 

The general idea of the proof is similar to that of the proof for $\text{Co}_1$. 
Suppose that $G$ is a finite group with $\Gamma(G) = \Gamma(\text{M})$. 
Then \cite[Theorem~3]{Hagie} implies that $G/F(G) \cong \text{M}$ and that $F(G)$ is either trivial or a $3$-group. 
We must show that $F(G)$ is trivial, so it suffices to show that $O_3(G)$ is trivial.  
Supposing towards a contradiction that $O_3(G)$ is non-trivial, we infer, by an argument completely analogous to the one given in Section~\ref{s:proofCo1}, that $\text{M}$ must admit a faithful irreducible module $V$ in characteristic $3$ on which every element of order $41$, $59$ or $71$ acts fixed-point freely. 
In particular, every element of order $41$ in a maximal subgroup $H = 3^8.\text{O}_8^-(3).2_3$ of $\text{M}$ must act fixed-point freely on this module $V$. 

We first claim that the only faithful irreducible modules for $H$ in characteristic $3$ on which elements of order $41$ act fixed-point freely have dimension $8$, $56$ or $104$, and that there is a unique module of each such dimension. 
The normal subgroup $3^8$ of $H$ coincides with $O_3(H)$, so \cite[Lemma 15.37]{Isaacs} implies that every irreducible $3$-modular representation of $H$ is unfaithful with image $H/O_3(H) \cong \text{O}_8^-(3).2_3$. 
Hence, the $3$-modular character table of $H$ is the same as that of $\text{O}_8^-(3).2_3$. 
That is to say, given a homomorphism $\varphi : H \rightarrow H/O_3(G)$ and an irreducible $3$-modular Brauer character $\rho$ for $H$, we have $\rho(x) = \rho'(\varphi(x))$ for all $x \in H$, where $\rho'$ is the character of the corresponding representation of $H/O_3(H)$.  
The $3$-modular character table of the index-$2$ subgroup $\text{O}_8^-(3)$ of $\text{O}_8^-(3).2_3$ is available in {\sf GAP}~\cite{GAPbc,GAP}, and we infer from it that elements of order $41$ in $\text{O}_8^-(3)$ act fixed-point freely only in the unique representations of dimensions $8$, $56$ and $104$. 
These representations extend to $\text{O}_8^-(3).2_3$ by uniqueness and Clifford's theorem, so the claim holds.

Next, we observe that the conjugacy class 2A in $\text{M}$ intersects $H$. 
One way to see this is to first recall that the only other involutions in $\text{M}$ have centraliser $2^{1+24}.\text{Co}_1$ (see \cite{ATLAS}). 
Constructing $H$ in {\sc Magma} \cite{Magma} using the generators given in the online ATLAS \cite{onlineATLAS}, we check directly that $H$ contains involutions centralised by a group of order $3^{12}$, whereas the highest power of $3$ dividing $|\text{Co}_1|$ is $3^9$. 
Now we check (using {\sc Magma} \cite{Magma}) that exactly four classes of elements of order $4$ in $H$ square to elements of the aforementioned class of involutions. 
Denote the union of these classes by $\mathcal{C}$, and note that each class in $\mathcal{C}$ belongs to the class 4B in $\text{M}$, because this is the only class of elements of order $4$ in $\text{M}$ that square to elements of 2A. 
Next, we define a homomorphism $\varphi : H \rightarrow H / O_3(H) \cong \text{O}_8^-(3).2_3$ using the \texttt{LMGRadicalQuotient} function in {\sc Magma}~\cite{Magma}, and check that there exist classes in $\mathcal{C}$ that project under $\varphi$ to the classes 4A and 4D in $\text{O}_8^-(3).2_3$. 
Let us fix labels 4A$'$ and 4D$'$ for these two classes, so that 4A$'$ projects to 4A and 4D$'$ projects to 4D. 
(Note that two of the classes in $\mathcal{C}$ project to 4A, but the following argument does not depend on which one we choose to work with. The final class projects to 4E.) 

Let $\rho$ denote the restriction to $H$ of the $3$-modular Brauer character corresponding to the hypothesised module $V$ for $\text{M}$. 
Then $\rho$ must be a linear combination with non-negative integer coefficients of the characters of the aforementioned $8$-, $56$- and $104$-dimensional modules for $H$. 
Label these characters and their associated coefficients by $\rho_d$ and $x_d$, respectively, for each $d \in \{8,56,104\}$. 
Since both of the classes 4A$'$ and 4D$'$ in $H$ belong to the same class (4B) in $\text{M}$, we must have $\rho (\text{4A}') = \rho(\text{4D}')$.
Since they project under the homomorphism $\varphi$ to the classes 4A and 4D in $\text{O}_8^-(3).2_3$, and (as explained above) the $3$-modular character tables of $H$ and $\text{O}_8^-(3).2_3$ coincide, we must have
\[
\rho_8'(4\text{A})x_8 + \rho_{56}'(4\text{A})x_{56} + \rho_{104}'(4\text{A})x_{104} = 
\rho_8'(4\text{D})x_8 + \rho_{56}'(4\text{D})x_{56} + \rho_{104}'(4\text{D})x_{104},
\]
where $\rho_d'$ denotes the $3$-modular Brauer character of $\text{O}_8^-(3).2_3$ corresponding (via $\varphi$) to $\rho_d$. 
Since the classes 4A and 4D in $\text{O}_8^-(3).2_3$ lie in $\text{O}_8^-(3)$, and (as explained above) the $8$-, $56$- and $104$-dimensional representations of $\text{O}_8^-(3)$ extend to $\text{O}_8^-(3).2_3$, we can read off the character values in the above equation from the $3$-modular character table of $\text{O}_8^-(3)$ (which is available in {\sf GAP}~\cite{GAPbc,GAP}) to deduce that
\[
-6x_8 - 26x_{56} - 38x_{104} = -2x_8 + 2x_{56} - 2x_{104}, 
\quad \text{i.e.} \quad 4x_8 + 28x_{56} + 36x_{104} = 0.
\]
Given that this equation has no (non-trivial) non-negative integer solutions $(x_8,x_{56},x_{104})$, the module $V$ does not exist. 
Therefore, $O_3(G)$ is trivial, and so $G \cong \text{M}$.

\section{Proof of Theorem~\ref{thm1} --- $\text{B}$} \label{s:proofB}

Recall that $|\text{B}| = 2^{41} \cdot 3^{13} \cdot 5^6 \cdot 7^2 \cdot 11 \cdot 13 \cdot 17 \cdot 19 \cdot 23 \cdot 31 \cdot 47$. 
We see from the ATLAS~\cite{ATLAS} that the prime graph $\Gamma(\text{B})$ has two isolated vertices, namely $31$ and $47$, and that all of its other vertices lie in a connected component. 
We also recall that from the ATLAS~\cite{ATLAS} that $\text{B}$ contains amongst its maximal subgroups the groups $\text{L}_2(31)$, $47:23$ and $2^{1+22}.\text{Co}_2$, and a group of shape $[2^{30}].\text{L}_5(2)$ (an extension of $\text{L}_5(2)$ by a certain group of order $2^{30}$.) 

The strategy here is similar to that of the preceding two sections, but the argument is a little more involved. 
Suppose that $G$ is a finite group with $\Gamma(G) = \Gamma(\text{B})$. 
Then \cite[Theorem~3]{Hagie} implies that $G/F(G) \cong \text{B}$ and that $F(G)$ is either trivial or a $2$-group. 
We must show that $F(G)$ is trivial, so it suffices to show that $O_2(G)$ is trivial. 
Supposing towards a contradiction that $O_2(G)$ is non-trivial, we infer that $\text{B}$ must admit a faithful irreducible module $V$ in characteristic $2$ on which every element of order $31$ or $47$ acts fixed-point freely. 
Let $\chi$ denote the corresponding Brauer character. 
We first consider restricting $V$ to the maximal subgroups $\text{L}_2(31)$ and $[2^{30}].\text{L}_5(2)$ of $\text{B}$, which contain elements of order $31$, and the maximal subgroup $47:23$, which contains elements of order $47$. 
This allows us to deduce the values of $\chi$ on various conjugacy classes that intersect the maximal subgroup $2^{1+22}.\text{Co}_2$. 
Since $2^{1+22}.\text{Co}_2$ has relatively few classes of odd-order elements, it turns out that we can impose sufficiently many constraints to establish the non-existence of $V$. 

We first consider $V \downarrow \text{L}_2(31)$. 
By the $2$-modular character table of $\text{L}_2(31)$, the only irreducible modules for $\text{L}_2(31)$ in characteristic $2$ on which elements of order $31$ act fixed-point freely are the two modules of dimension $15$. 
Let $\chi_2$ and $\chi_3$ denote the corresponding Brauer characters, for consistency with the numbering in the character table as given in {\sf GAP}~\cite{GAPbc,GAP}. 
Then there exist non-negative integers $k_2$ and $k_3$ (not both zero) such that
\begin{equation} \label{chiB31}
\chi |_{\text{L}_2(31)} = k_2\chi_2 + k_3\chi_3.
\end{equation}
Moreover, $\chi_2$ and $\chi_3$ vanish on all elements of order $3$, $5$ and $15$ in $\text{L}_2(31)$, and all such elements belong to the conjugacy classes 3B, 5B and 15B in $\text{B}$, respectively, so 
\begin{equation} \label{3B-5B-15B}
\chi(3\text{B}) = \chi(5\text{B}) = \chi(15\text{B}) = 0.
\end{equation}

Next, we claim that $k_2 = k_3$. 
To prove this, consider one of the two conjugacy classes of elements of order $31$ in $\text{L}_2(31)$, say the class 31A, which belongs to the class 31A in $\text{B}$.  
Checking the values of $\chi_2$ and $\chi_3$ on 31A elements in $\text{L}_2(31)$, we infer that the value of $\chi$ on 31A elements in $\text{B}$ must be
\begin{equation} \label{x31A}
\chi(31\text{A}) = k_2c + k_3\overline{c} \quad \text{where} \quad c = \frac{-1+\sqrt{-31}}{2}
\end{equation}
and $\overline{c}$ denotes the complex conjugate of $c$. 
Now consider $V \downarrow H$, where $H = [2^{30}].L_5(2)$. 
There are $16$ irreducible $2$-modular Brauer characters for $H$. 
Label them $\rho_1,\dots,\rho_{16}$ in accordance with the numbering in the $2$-modular character table of $H$ in {\sf GAP}~\cite{GAPbc,GAP}, and let $x_1,\dots,x_{16}$ denote the corresponding multiplicities in the restriction of $\chi$ to $H$. 
The requirement that all elements of order $31$ in $H$ act fixed-point freely on $V$ implies that $x_i=0$ for all $i \in \{1,6,11,12,13,14,15,16\}$. 
In particular, we have
\begin{equation} \label{x31A-2}
\chi |_H = \sum_{i \in I} x_i \rho_i \quad \text{where} \quad I = \{2,3,4,5,7,8,9,10\}.
\end{equation}
Note that the remaining possibilities for the composition factors of $V \downarrow H$ are the two modules of dimension $5$, corresponding to $i \in \{2,3\}$, the two modules of dimension $10$, corresponding to $i \in \{4,5\}$, and the four modules of dimension $40$, corresponding to $i \in \{7,8,9,10\}$. 
Now, consider the class 31A in $H$. 
Elements in this class also belong to the class 31A in $\text{B}$, so we can use \eqref{x31A-2} to obtain an expression for $\chi(31A)$ in terms of the known values of the characters $\rho_i$, $i \in I$, on 31A elements of $H$, and the unknown multiplicities $x_i$, $i \in I$. 
The values of the $\rho_i$ on 31A elements of $H$ are all integer linear combinations of primitive 31st roots of unity, and so too are the complex numbers $c$ and $\overline{c}$ in \eqref{x31A}. 
Hence, equating the right-hand sides of \eqref{x31A} and \eqref{x31A-2} (with the latter applied to 31A elements of $H$), we obtain an equality between two integer linear combinations of primitive 31st roots of unity; one involving the unknowns $k_2$ and $k_3$, the other involving the unknowns $x_i$, $i \in I$. 
Since the primitive 31st roots of unity are linearly independent over $\mathbb{Z}$, we can regard this equality as a system of $30$ linear equations (one per root) in the eight unknowns $x_i$, $i \in I$. 
It turns out that only six of these equations are distinct; explicitly, we obtain the linear system with augmented matrix
\[
\left[
\begin{array}{rrrrrrrr}
1&\cdot&\cdot&\cdot&1&\cdot&1&2 \\
\cdot&\cdot&1&\cdot&1&2&1&1 \\
\cdot&\cdot&1&\cdot&2&2&1&2 \\
\cdot&\cdot&\cdot&1&2&1&1&1 \\
\cdot&\cdot&\cdot&1&2&2&2&1 \\
\cdot&1&\cdot&\cdot&\cdot&1&2&1 \\
\end{array}
\right| \hspace{-2pt}
\left.
\begin{array}{r}
k_2 \\ k_3 \\ k_2 \\ k_2 \\ k_3 \\ k_3
\end{array} 
\right]
\sim
 \left[
\begin{array}{rrrrrrrr}
1&\cdot&\cdot&\cdot&\cdot&\cdot&1&1 \\
\cdot&1&\cdot&\cdot&\cdot&\cdot&1&1 \\
\cdot&\cdot&1&\cdot&\cdot&\cdot&-1&\cdot \\
\cdot&\cdot&\cdot&1&\cdot&\cdot&\cdot&-1 \\
\cdot&\cdot&\cdot&\cdot&1&\cdot&\cdot&1 \\
\cdot&\cdot&\cdot&\cdot&\cdot&1&1&\cdot \\
\end{array}
\right| \hspace{-2pt}
\left.
\begin{array}{c}
k_3 \\ k_2 \\ k_2 \\ k_3 \\ k_2-k_3 \\ k_3-k_2
\end{array} 
\right],
\]
where $\cdot=0$ and $\sim$ denotes row equivalence. 
Since the $x_i$ must be non-negative integers, it follows that 
\begin{equation} \label{31sol}
x_7=x_8=x_9=x_{10}=0 \quad \text{and} \quad x_2=x_3=x_4=x_5=k_2=k_3.
\end{equation}
In particular, $k_2=k_3$ as claimed. 

Let us therefore write $k=k_2=k_3$. 
Note first that \eqref{chiB31} then reads $\chi |_{\text{L}_2(31)} = k(\chi_2 + \chi_3)$, and so the value of $\chi$ on the identity class 1A of $\text{B}$ is
\begin{equation} \label{1A}
\chi(1A) = 30k, 
\end{equation}
because (as noted above) the modules for $\text{L}_2(31)$ corresponding to the characters $\chi_2$ and $\chi_3$ both have dimension $15$. 
Now apply \eqref{31sol} also to \eqref{x31A-2}. 
This yields
\begin{equation} \label{finalChiH}
\chi |_H = k(\rho_2+\rho_3+\rho_4+\rho_5).
\end{equation}
In other words, $V$ restricts to $H$ as $k$ copies of each of the two $5$-dimensional and two $10$-dimensional modules for $H$. 
Next, consider either of the classes 7A and 7B in $H$, both of which belong to the class 7A in $\text{B}$. 
By \eqref{finalChiH} and the $2$-modular character table of $H$, 
\begin{equation} \label{7A}
\chi(7A) = k \left( \frac{3+\sqrt{-7}}{2} + \frac{3-\sqrt{-7}}{2} + \frac{-1+\sqrt{-7}}{2} + \frac{-1-\sqrt{-7}}{2} \right) = 2k. 
\end{equation}
Similarly, the class 3A in $H$ belongs to the class 3A in $\text{B}$, so 
\begin{equation} \label{3A}
\chi(3A) = k (2+2+1+1) = 6k.
\end{equation}
Finally, we claim that
\begin{equation} \label{23AB}
\chi(23A) = \chi(23B) = 0.
\end{equation}
This is easily seen by considering the restriction of $V$ to $47:23$. 
Since elements of order $47$ in $\text{B}$ must act fixed-point freely on $V$, the only possible composition factors of $V \downarrow 47:23$ are the two modules of dimension $23$, and the Brauer characters of these modules vanish on all elements of order $23$ in $47:23$. 
In particular, $\chi$ must satisfy \eqref{23AB}. 

\begin{table}[!t]
\begin{tabular}{llr}
\toprule
$\text{B}$ & $2^{1+22}.\text{Co}_2$ & $\chi$ \\
\midrule
1A & 1A & $30k$ \\
3A & 3B & $6k$ \\
3B & 3A & $0$ \\
5B & 5A & $0$ \\
7A & 7A & $2k$ \\
15B & 15B, 15C & $0$ \\
23A & 23A & $0$ \\
23B & 23B & $0$ \\
\bottomrule
\end{tabular}
\caption{Some conjugacy classes of odd-order elements in $\text{B}$ that intersect the maximal subgroup $2^{1+22}.\text{Co}_2$, the corresponding classes in $2^{1+22}.\text{Co}_2$, and the corresponding values of the $2$-modular Brauer character $\chi$ of the hypothesised module $V$ for $\text{B}$ in Section~\ref{s:proofB}, per \eqref{3B-5B-15B}, \eqref{1A}, \eqref{7A}, \eqref{3A} and \eqref{23AB}.}
\label{tab2}
\end{table}

We now complete the proof by considering $V \downarrow K$, where $K = 2^{1+22}.\text{Co}_2$. 
Table~\ref{tab2} summarises the supposed values of $\chi$ on various conjugacy classes of $\text{B}$, per \eqref{3B-5B-15B}, \eqref{1A}, \eqref{7A}, \eqref{3A} and \eqref{23AB}, and lists the corresponding classes in $K$. 
Now, $K$ has only $13$ irreducible $2$-modular Brauer characters. 
Denote the associated coefficients in $\chi |_K$ by $x_1,\dots,x_{13}$, following the numbering in the $2$-modular character table of $K$ in {\sf GAP}~\cite{GAPbc,GAP}. 
Using that character table and Table~\ref{tab2}, we obtain a system of nine equations (one per $K$-class appearing in Table~\ref{tab2}) in the $13$ unknowns $x_1,\dots,x_{13}$, with augmented matrix
\[
\small
\hspace{-2pt}
\left[
\begin{array}{rrrrrrrrrrrrr}
1&22&230&748&748&3520&5312&8602&8602&36938&83948&156538&1835008 \\
1&-5&14&-8&-8&-44&20&43&43&-79&-22&100&-128 \\
1&4&5&1&1&10&2&-20&-20&-52&-58&-80&-128 \\
1&-3&5&-2&-2&-5&12&2&2&13&-2&-12&8 \\
1&1&-1&-1&-1&-1&-1&-1&-1&-1&-3&-3&\cdot \\
1&\cdot&-1&a&\overline{a}&1&\cdot&-a&-\overline{a}&1&-2&\cdot&2 \\
1&\cdot&-1&\overline{a}&a&1&\cdot&-\overline{a}&-a&1&-2&\cdot&2 \\
1&-1&\cdot&b&\overline{b}&1&-1&\cdot&\cdot&\cdot&-2&\cdot&-1 \\
1&-1&\cdot&\overline{b}&b&1&-1&\cdot&\cdot&\cdot&-2&\cdot&-1
\end{array}
\right| \hspace{-3pt}
\left.
\begin{array}{r}
30k \\
\cdot \\
6k \\
\cdot \\
2k \\
\cdot \\
\cdot \\
\cdot \\
\cdot
\end{array} 
\right]
\]
where $\cdot = 0$ and
\[
a = \frac{-1+\sqrt{-15}}{2}, \quad b = \frac{1-\sqrt{-23}}{2}.
\]
This system has no non-negative integer solutions $(x_1,\dots,x_{13})$ for any positive integer $k$ (see below), so the hypothesised module $V$ for $\text{B}$ does not exist. 
Therefore, $O_2(G)$ is trivial, so $F(G)$ is trivial and hence $G \cong \text{B}$, as claimed. 

It remains to verify that the aforementioned linear system has no non-negative integer solutions. 
It suffices to consider $k=1$ and check that there are no non-negative rational solutions. 
When $k=1$, the augmented matrix given above is row equivalent to
\[
\left[
\begin{array}{rrrrrrrrrrrrr}
1&\cdot&\cdot&\cdot&\cdot&\cdot&\cdot&\cdot&\cdot&12&20&34&468 \\
\cdot&1&\cdot&\cdot&\cdot&\cdot&\cdot&\cdot&\cdot&7&20&36&478 \\
\cdot&\cdot&1&\cdot&\cdot&\cdot&\cdot&\cdot&\cdot&-26&-30&-34&-446 \\
\cdot&\cdot&\cdot&1&\cdot&\cdot&\cdot&\cdot&\cdot&23&35&49&645 \\
\cdot&\cdot&\cdot&\cdot&1&\cdot&\cdot&\cdot&\cdot&23&35&49&645 \\
\cdot&\cdot&\cdot&\cdot&\cdot&1&\cdot&\cdot&\cdot&-12&-16&-20&-260 \\
\cdot&\cdot&\cdot&\cdot&\cdot&\cdot&1&\cdot&\cdot&16&21&27&376 \\
\cdot&\cdot&\cdot&\cdot&\cdot&\cdot&\cdot&1&\cdot&-2&-1&1&-7 \\
\cdot&\cdot&\cdot&\cdot&\cdot&\cdot&\cdot&\cdot&1&-2&-1&1&-7
\end{array}
\right| \hspace{-2pt}
\left.
\begin{array}{r}
19982/9315 \\
17422/9315 \\
-3668/3105 \\
1982/1035 \\
1982/1035 \\
-11176/9315 \\
3074/3105 \\
-1972/9315 \\
-1972/9315
\end{array} 
\right].
\]
Given that we need $x_7 \geqslant 0$, the third-last row forces 
$x_{10} < 1/16$, $x_{11} < 1/21$ and $x_{13} < 1/376$. 
Putting this and $x_{12} \geqslant 0$ into the last row yields $x_9 < 0$.

\section*{Acknowledgements} 
This work was initiated in preparation for a research retreat run by the The University of Auckland's Algebra \& Combinatorics research group, on Waiheke Island in July 2021. 
We thank all of our fellow participants for an enjoyable and productive retreat, and in particular for discussions with various colleagues about the recognisability-by-prime-graph problem. 
Special thanks are due to Jeroen~Schillewaert for organising the retreat, and to Eamonn O'Brien and Gabriel Verret for helpful discussions about the questions that we have managed to answer here. 
We are also grateful to Chris Parker and an anonymous referee for several suggestions that helped to improve the paper.

\end{document}